\documentclass{amsart}
\usepackage{amssymb,latexsym}
\theoremstyle{plain}
\newtheorem{theorem}{Theorem}

\newtheorem{lemma}{Lemma}

\theoremstyle{definition}

\theoremstyle{remark}

\numberwithin{equation}{section}
\newcommand{\s}{\sigma}
\newcommand{\R}{\mathbb R}
\newcommand{\Rn}{\mathbb R^n}
\newcommand{\RM}{\mathbb R^{n+1}_+}
\newcommand{\Rm}{\mathbb R^{n+1}}
\newcommand{\al}{\alpha}
\newcommand{\e}{\epsilon}
\begin{document}
\title[Decay at infinity of caloric functions\dots]{Decay at infinity of caloric
functions\\within characteristic hyperplanes}
\author{L. Escauriaza}
\address[L. Escauriaza]{Universidad del Pa{\'\i}s Vasco / Euskal Herriko
Unibertsitatea\\Dpto. de Matem\'aticas\\Apto. 644, 48080 Bilbao, Spain.} 
\email{mtpeszul@lg.ehu.es}
\thanks{The first and fourth authors are supported  by MEC grant MTM2004-03029 and by the
European Commission via the network Harmonic Analysis and Related Problems, project number
RTN2-2001-00315. The second and third authors are supported by NFS grants. }
\author{C.E. Kenig}
\address[C.E. Kenig]{Department of Mathematics\\University of Chicago\\5734 S. University Avenue\\Chicago,
Illinois 60637, USA.} 
\email{cek@math.uchicago.edu}
\author{G. Ponce}
\address[G. Ponce]{Department of Mathematics\\ South Hall, Room 6607\\ University of California, Santa Barbara\\CA 93106, USA.} 
\email{ponce@math.ucsb.edu}
\author{L. Vega}
\address[L. Vega]{Universidad del Pa{\'\i}s Vasco / Euskal Herriko
Unibertsitatea\\Dpto. De matem\'aticas\\Apto. 644, 48080 Bilbao, Spain.} 
\email{mtpvegol@lg.ehu.es}
\begin{abstract} It is shown that a function $u$ satisfying,
$|\Delta u+\partial_tu|\le M\left(|u|+|\nabla u|\right)$,
$|u(x,t)|\le Me^{M|x|^2}$ in
$\R^n\times [0,T]$ and $|u(x,0)|\le C_ke^{-k|x|^2}$ in $\R^n$ and for all  
 $k\ge 1$, must vanish identically in
$\R^n\times [0,T]$.
\end{abstract}
\maketitle

\section{Introduction}\label{S:1}
E.M. Landis and O.A. Oleinik asked \cite[\S4]{lo} for a proof of the
following conjecture:

\medskip
{\it If $u(x,t)$ is a bounded solution of a uniformly parabolic equation
\begin{equation*}
Pu=\sum_{i,j=1}^n\partial_i\left(g^{ij}(x)\partial_ju\right)-\partial_tu+b(x)\cdot\nabla
u+c(x)u=0\ ,
\end{equation*}
in the layer $\R^n\times [0,T]$ and the condition, $|u(x,T)|\le
Ne^{-|x|^{2+\e}}$, $x\in\R^n$, holds for some positive constants $N$ and $\e$, then
$u(x,t)\equiv 0$ in $\R^n\times [0,T]$.}
\medskip

As they wrote it, natural conditions should be placed on the behavior of the coefficients of $P$
at infinity for the conjecture to hold. 

Here, we give an answer to this question when the leading parabolic operator is the
backward heat operator and the lower order coefficients are bounded. In particular,
we prove the following quantitative and qualitative results of unique
continuation:   
\begin{theorem}\label{T:1} 
Assume that a function $u$ verifies the inequalities 
\begin{equation}\label{E:1}
|\Delta u+\partial_tu|\le M\left(|u|+|\nabla u|\right)\ \text{and}\ \  |u(x,t)|\le
Me^{M|x|^2}\ \text{in}\ \ \R^n\times [0,T]\ .
\end{equation}
Then, the following holds:
\begin{itemize}
\item If $\|u(\,\cdot\,,0)\|_{L^2\left(B_1\right)}$ is positive, there is
$N>0$ such that, when $|y|\ge N$
\begin{equation*}
\|u(\,\cdot\,,0)\|_{L^2\left(B_{|y|/2}(y)\right)}\ge
e^{-N|y|^2}\ \text{and}\ \|u(\,\cdot\,,0)\|_{L^2\left(B_1(y)\right)}\ge
e^{-N|y|^2\log{|y|}}\ . 
\end{equation*}
\item $u\equiv 0$ in $\R^n\times [0,T]$ if $|u(x,0)|\le C_ke^{-k|x|^2}$ for all $x\in\R^n$
and $k\ge 1$.
\end{itemize}  
\end{theorem}
Here, $B_r(y)=\{x\in\R^n:|x-y|< r\}$ and $B_r=B_r(0)$. We work with backward parabolic
operators because it is more convenient in this context.

When $n\ge 2$, we understand how to establish the conjecture when the matrix of
coefficients of the parabolic operator
$P$ verifies, for some large $M>0$, the conditions
\[
M^{-1}|\xi|^2\le\sum_{i,j=1}^ng^{ij}(x)\xi_i\xi_j\le M|\xi|^2\ \text{and}\ e^{M|x|}|\nabla
g^{ij}(x)|\le M\ ,\ \text{when}\ x,\ \xi\in\R^n\ ,
\]
and this will appear
elsewhere. When $n=1$ and provided that, $M^{-1}\le \gamma(x)\le M$ and $|\gamma'(x)|\le M$,
the changes of variables
\[y=\int_0^x\frac{ds}{\sqrt{\gamma(s)}}\ ,\ v(x,t)=u(y,t)\]
transform solutions $v$ of the
inequalities
\[|\partial_x\left(\gamma(x)\partial_xv\right)+\partial_tv|\le M\left(|v|+|\partial_x v|\right)\ ,\ |v(x,t)|\le
Me^{M|x|^2}\ \text{in}\ \ \R\times [0,T]\ ,\]
into solutions $u$ of backward parabolic inequalities, where the leading  operator
is the backward heat operator, as in \eqref{E:1}. This and Theorem
\ref{T:1} prove the conjecture when $n=1$.

The first author,
G. Seregin and V.
\v Sver\'ak proved in
\cite{ess2} the following qualitative property of unique continuation:

\medskip
{\it Let $\R^n_+=\{x=(x',x_n)\in \R^n : x_n>0\}$ and assume that $u$ satisfies 
\begin{equation}\label{E:2}
|\Delta u+\partial_tu|\le M\left(|u|+|\nabla u|\right)\ \text{,}\ \  |u(x,t)|\le
Me^{M|x|^2}\ \text{in}\ \ \R^n_+\times [0,T]\
\end{equation}
and $u(x,0)=0$ in $\R^n_+$. Then, $u\equiv 0$ in $\R^n_+\times [0,T]$}.
\medskip

This result is of interest in control theory; see \cite{mz}, and as explained in \cite{ssv} and \cite{ess3},
results of this type have shown to be helpful in the regularity theory for the Navier-Stokes
equations. The arguments in the proof of Theorem \ref{T:1} also imply the
following improvement of the last result.

\begin{theorem}\label{T:2} 
Let $u$ verify \eqref{E:2} and set $e_n=(0,\,\dots,0,1)$. Then, the following holds:
\begin{itemize}
\item If $\|u(\,\cdot\,,0)\|_{L^2\left(B_1(4e_n)\right)}$ is positive, there is
$N>0$ such that, when $y\ge N$
\begin{equation*}
\|u(\,\cdot\,,0)\|_{L^2\left(B_{y/2}(ye_n)\right)}\ge 
e^{-Ny^2}\ \text{and}\ \|u(\,\cdot\,,0)\|_{L^2\left(B_1(ye_n)\right)}\ge e^{-Ny^2\log{y}}\ .
\end{equation*}
\item $u\equiv 0$ in $\R^n_+\times [0,T]$ if $|u(x,0)|\le C_ke^{-k|x|^2}$ for all
$x\in\R^n_+$ and $k\ge 1$. 
\end{itemize}  
\end{theorem}

We present in sections \ref{S:2} and \ref{S:3} two different proofs of Theorems
\ref{T:1}. The first one is based on Carleman inequality
methods while the second on frequency function methods. The main tools in both proofs are a rescaling argument
and a quantification of the size of the constants involved in the two sphere and one cylinder inequalities
satisfied by solutions of certain parabolic equations, in terms of the $L^\infty$-norm of the lower order
coefficients and of the time of existence of solutions. See \cite[Lemma 3.10]{bk}, where similar ideas
appeared but dealing with three sphere inequalities and elliptic equations. In section
\ref{S:4}, we outline the proof of Theorem \ref{T:2}.

With the purpose of simplifying
 the arguments below, we only prove Theorems \ref{T:1} and \ref{T:2},
when the growth condition  in \eqref{E:1} or \eqref{E:2}, 
$|u(x,t)|\le Me^{M|x|^2}$,
is replaced by $u$ is bounded. The interested reader can easily
verify that the arguments below can be adapted to the more general case. 
\section{First Proof of Theorem \ref{T:1}}\label{S:2} 
The next five Lemmas are used in the first proof of Theorem \ref{T:1}.
The first one, Lemma \ref{L:decay}, is in a certain sense a localized version of the standard
energy inequality satisfied by solutions of parabolic inequalities (See \cite[Lemmas 1 and
5]{efv} for other versions of this Lemma). The Lemmas \ref{L:ineq} and
\ref{L:implidoub} appeared in \cite[Lemmas 2, 3]{efv}. 

\begin{lemma}\label{L:decay}
Assume that $u$ satisfies, $|\Delta u+\partial_tu|\le  R^2|u|+R|\nabla u|$, 
$\|u\|_{\infty}\le 1$,
$\|\nabla u\|_\infty\le R$ in $B_4\times [0,\frac 1{R^2}]$ and 
$\|u(\,\cdot\,,0)\|_{L^2\left(B_\rho\right)}\ge \theta R^{-n/2}$
for some $\theta$, $\rho$ in $(0,1]$ and $R>0$. Then, there is $N=N(n,\theta)$ such that
the  inequality, $\sqrt{N}\|u(\,\cdot\,,t)\|_{L^2\left(B_{2\rho}\right)}\ge R^{-n/2}$ holds,
when
$0<t\le 1/R^2$ and $R> N/\rho$. 

\end{lemma}
\begin{proof}
 Assume first that $\rho=1$ and set $f=u{\varphi}$, where
$\varphi\in C_0^\infty(B_2)$,
$0\le\varphi
\le 1$ and
$\varphi=1$ in
$B_{3/2}$. Then,
\begin{equation}\label{E:heat2}
|\Delta f+\partial_tf|\le R^2|f|+R|\nabla f|+NR\chi_{B_2\setminus B_{3/2}}\ .
\end{equation}
Setting $H(t)=\int f^2(x,t)G(x-y,t)\,dx$,
where
$G(x,t)=t^{-n/2}e^{-|x|^2/4t}$ and
$y\in B_1$, we have
\begin{equation}\label{E:heat3}\dot H(t)=2\int f(\Delta
f+\partial_tf)G(x-y,t)\,dx+2\int |\nabla f|^2G(x-y,t)\,dx\ ,\end{equation}
and from \eqref{E:heat2}, \eqref{E:heat3} and the Cauchy-Schwarz's inequality
\begin{equation*}\dot H(t)\ge -8R^2H(t)-Ne^{-1/Nt}\quad .\end{equation*}
Integration of this inequality in $(0,t)$, $0< t\le \frac 1{R^2}$,  gives
\begin{equation*}
N\int f^2(x,t)G(x-y,t)\,dx\ge u^2(y,0)-Ne^{-1/Nt}\quad .
\end{equation*}
Integrating the last inequality over $B_1$ and recalling that $\int
G(x-y,t)\,dy=1$, we get
\begin{equation*}
N\int_{B_2}u^2(x,t)\,dx\ge \int_{B_1}u^2(x,0)\,dx-Ne^{-1/Nt}\ge
R^{-n}\left(\theta^2-Ne^{-R^2/2N}\right)\ ,
\end{equation*}
when $0<t\le 1/R^2$, which implies Lemma \ref{L:decay} when $\rho=1$. 

When $\rho$ is in $(0,1)$,  the function $u_\rho(x,t)=u(\rho x,\rho^2t)$,  satisfies the conditions  in 
Lemma \ref{L:decay} with
$\rho=1$
 and $R$ replaced by $\rho R$. The Lemma then, follows after rescaling to the case $\rho=1$.
\end{proof}
\begin{lemma}\label{L:ineq} The inequality
\begin{equation*}
\int \tfrac {|x|^2}{8a}h^2e^{-|x|^2/4a}\,dx\le 2a\int |\nabla
h|^2e^{-|x|^2/4a}\,dx+\tfrac n2\int h^2e^{-|x|^2/4a}\,dx
\end{equation*}
holds for all $h\in C_0^\infty(\Rn)$ and $a>0$.
\end{lemma}
 \begin{proof}
 The inequality follows setting  $v=he^{-|x|^2/8a}$ and from the identity
\begin{align*}
2a\int |\nabla
h|^2e^{-|x|^2/4a}\,dx+&\tfrac n2\int h^2e^{-|x|^2/4a}\,dx-\int \tfrac
{|x|^2}{8a}h^2e^{-|x|^2/4a}\,dx\\&=2a\int |\nabla v|^2\,dx\ .
\end{align*}
\end{proof}
\begin{lemma}\label{L:implidoub} Assume that $N$ and $\Theta$ verify
$N\log(N\Theta)\ge 1$,
$h\in C_0^\infty(\Rn)$ and that the inequality 
\begin{equation*}
\begin{split}
2a\int |\nabla h|^2e^{-|x|^2/4a}\,dx&+\tfrac n2\int
h^2e^{-|x|^2/4a}\,dx\le N\log{(N\Theta)}\int
h^2e^{-|x|^2/4a}\,dx
\end{split}
\end{equation*}
holds, when $0<a\le\frac 1{12N\log{(N\Theta)}}$. Then,
\begin{equation*}
\int_{B_{2r}}h^2\,dx\le (N\Theta)^N\int_{B_r}h^2\,dx\ ,\ \text{when}\ 0<r\le 1/2\ .
\end{equation*}
\end{lemma}
\begin{proof}
The inequality satisfied by $h$ and Lemma \ref{L:ineq} show that  
\begin{equation*}
\int \tfrac{|x|^2}{8a}h^2e^{-|x|^2/4a}\,dx\le N\log{(N\Theta)}\int
h^2e^{-|x|^2/4a}\,dx
\end{equation*}
when $a\le1/\left(12N\log{(N\Theta)}\right)$.
For given $0<r\le 1/2$ and $0< a\le
\tfrac{r^2}{16N\log{\left(N\Theta\right)}}\
$, the last inequality  implies 
\begin{gather*}
\int\tfrac{|x|^2}{8a}h^2e^{-|x|^2/4a}\,dx\le N\log{(N\Theta)}\left[\int_{B_r}
h^2\,dx+\tfrac{8a}{r^2}\int_{\Rn\setminus B_r}
\tfrac{|x|^2}{8a}h^2e^{-|x|^2/4a}\,dx\right]\\\le N\log{(N\Theta)}\int_{B_r} h^2\,dx+\tfrac
12\int
\tfrac{|x|^2}{8a}h^2e^{-|x|^2/4a}\,dx\ .
\end{gather*}
Thus,
\begin{equation}\label{E:casiadob}
\int \tfrac{|x|^2}{16a}h^2e^{-|x|^2/4a}\,dx\le N\log{(N\Theta)}\int_{B_r}
h^2\,dx\ ,\  \text{when}\  0< a\le
\tfrac{r^2}{16N\log{\left(N\Theta\right)}}\ .
\end{equation}
Now, $e^{-|x|^2/4a}|x|^2/(16a)\ge
(N\Theta)^{-N}N\log{(N\Theta)}$ when $r\le |x|\le 2r$ and $a=
\tfrac{r^2}{16N\log{\left(N\Theta\right)}}\ $. This and \eqref{E:casiadob}
imply
\begin{equation*}
\int_{B_{2r}}h^2\,dx\le (N\Theta)^N\int_{B_r}h^2\,dx\ ,\ \text{when}\ 0<r\le 1/2\ .
\end{equation*}
\end{proof}

The Lemma \ref{L:carleman1} contains the Carleman inequality we need. Here, $dX=dxdt$ is the
Lebesgue measure in
$\R^{n+1}_+$ and $\s_a(t)=\s(t+a)$, denotes the translation by $a>0$ of a function $\s$ of the
time-variable.
\begin{lemma}\label{L:carleman1}  Given $\al\ge 2+n/2$, there are $N=N(n)$ and an increasing function,
$\sigma :[0,+\infty)\longrightarrow [0,+\infty)$ verifying,
$t/N\le\sigma(t)\le t$ in
$[0,4/\alpha]$ and such that the inequality 
\begin{gather*}
\begin{split}
\al^2\int\s_a^{-\al}f^2e^{-|x|^2/4(t+a)}\,dX+\al\int \s_a^{1-\al}\ |\nabla
f|^2e^{-|x|^2/4(t+a)}\,dX&\\\le N\int\s_a^{1-\al}\left(\Delta
f+\partial_tf\right)^2e^{-|x|^2/4(t+a)}\,dX&
\end{split}\\ +\ 
\s(a)^{-\al}\left[-(a/N)\int|\nabla
f(x,0)|^2e^{-|x|^2/4a}\,dx+N\al\int f^2(x,0)e^{-|x|^2/4a}\,dx\right]
\end{gather*}
holds, when $0 < a\le\frac 1{\al}$ and $f\in C_0^\infty(\Rn\times[0,\frac 4\al))$. 
\end{lemma}

This inequality appeared
first in \cite[\S3]{f} in the context of variable
coefficients parabolic operators. The inequality is not stated there as it is shown
above, some additional terms appear or are missing on the right hand side of the corresponding
inequality in
\cite[\S3]{f}. These additional  terms arise from the purpose of controlling certain error terms generated by
the variable coefficients of the parabolic operator, but  they can be dropped when
the operator is the backward heat operator. Other versions of this inequality
appeared in
\cite[(1.4)]{ess1}, \cite[Proposition 6.1]{ess3} and \cite[\S3]{ess2} but none of them is
stated or proved as we need need it here.

As it is usual in the context of $L^2$-Carleman estimates,
we use suitable integration by parts to prove Lemma \ref{L:carleman1}. The calculations can be organized
either by using identities developed in
\cite[Lemma 1]{ef} and \cite[Lemma 3]{ess1}, or by following more or less standard calculations
with new dependent variables and commutators  in the spirit of \cite{h}, \cite{h1} or \cite{t}.
In this paper we will use the former method.
\begin{proof}  Assume first that the following claim holds: 

There are $N=N(n)$ and an increasing function,
$h:[0,+\infty)\longrightarrow [0,+\infty)$, verifying, $t/N\le h(t)\le t$ in $[0,6]$ and such
that the inequality
\begin{multline}\label{E:primera desigualdad}
\al\int h_a^{-\al}u^2e^{-|x|^2/4(t+a)}\,dX+\int h_a^{1-\al}\ |\nabla
u|^2e^{-|x|^2/4(t+a)}\,dX\\\le N\int h_a^{1-\al}\left(\Delta
u+\partial_tu\right)^2e^{-|x|^2/4(t+a)}\,dX\\ +\ 
h(a)^{-\al}\left[-(a/N)\int|\nabla
u(x,0)|^2e^{-|x|^2/4a}\,dx+N\al\int u^2(x,0)e^{-|x|^2/4a}\,dx\right]
\end{multline}
holds, when $\al\ge 2+n/2$, $0 < a\le 1$ and $u\in C_0^\infty(\Rn\times[0, 4))$.

Take as $u$ in \eqref{E:primera desigualdad} the function, $u(x,t)=f(x/\sqrt\al
,t/\al)$, when $f\in C_0^\infty(\Rn\times[0,\frac 4\al))$ and  define
$\s(t)=h(\al t)/\al$. Then,  it is simple to verify
that Lemma \ref{L:carleman1} holds after undoing the change of variables and counting of
the number of $\al$'s at each side of the inequality.

  In order to prove the claim we recall the following identity \cite[(2.4)]{ess1}, which
holds when $\al\in\R$,
$u\in C_0^\infty(\Rn\times[0, 4))$, $G$ is a positive caloric function in $\R^{n+1}_+$ and
$\gamma :[0,+\infty)\longrightarrow (0,+\infty)$ is an increasing
 smooth function: 
\begin{multline}\label{E:identidad}
\tfrac{2\gamma^{1-\al}}{\dot \gamma}\left(\partial_tu-\nabla\log G\cdot\nabla
u-\tfrac{\al\dot \gamma}{2\gamma}u\right)^2G+\tfrac{\gamma^{1-\al}}{\dot
\gamma}\mathcal{D}_G\nabla u\cdot\nabla u\ G\\=\tfrac{2\gamma^{1-\al}}{\dot
\gamma}\left(\partial_tu-\nabla\log G\cdot\nabla u-\tfrac{\al\dot
\gamma}{2\gamma}u\right)\left(\Delta
u+\partial_tu\right)G\\+\partial_t\left[\tfrac{\gamma^{1-\al}}{\dot \gamma}|\nabla
u|^2G-\tfrac{\al \gamma^{-\al}}{2}u^2G\right]\\ +\tfrac{\gamma^{1-\al}}{\dot
\gamma}\nabla\cdot\left[2\partial_tuG\nabla u+|\nabla u|^2\nabla G-2\left(\nabla G\cdot\nabla
u\right)\nabla u-\tfrac{\al\dot \gamma}{\gamma}uG\nabla u-\tfrac{\al\dot \gamma}{2\gamma}
u^2\nabla G\right]\  .
\end{multline}
Here, $\mathcal{D}_G$ denotes the $n\times n$ matrix
\begin{equation*}
\mathcal{D}_G=\overbrace{\log{\left(\tfrac{\gamma}{\dot \gamma}\right)}}^.\ \mathcal{I}+2D^2(\log G)\ .
\end{equation*}

If in \eqref{E:identidad} we set $\gamma(t)=h_a(t)$, where $h(t)=te^{-t/6}$, $a\in
(0,1]$ and let $G$ be the translated Gauss Kernel,
$G_a(x,t)=(t+a)^{-n/2}e^{-|x|^2/4(t+a)}$, we have
\begin{equation}\label{E:acotaciones}\tfrac 1{e}(t+a)\le h_a(t)\le t+a\ ,\   \tfrac
1{6e}\le\dot h_a(t)\le 1\ \text{and}\ \mathcal{D}_{G_a}\ge
\tfrac 16\mathcal{I},\ \text{when}\ t\in (0,4]\ .
\end{equation}
Integrating the
identity \eqref{E:identidad} over
$\Rm_+$, one gets from \eqref{E:acotaciones} and the Cauchy-Schwarz inequality (which is used  
to handle the first integral on the right hand side of the formula \eqref{E:identidad}), the 
 bound
\begin{multline}\label{E:carleman2}
\int h_a^{1-\al}|\nabla u|^2G_a\,dX\le N\int h_a^{1-\al}\left(\Delta
u+\partial_tu\right)^2G_a\,dX\\+\ 
h(a)^{-\al-n/2}\left[-(a/N)\int|\nabla
u(x,0)|^2e^{-|x|^2/4a}\,dx+N\al\int u^2(x,0)e^{-|x|^2/4a}\,dx\right]\  .
\end{multline}
Finally, the claim follows after multiplication of the identity
$$(\Delta +\partial_t)(u^2)=2u(\Delta u +\partial_t u)+2|\nabla u|^2$$
by $h_a^{1-\alpha}G_a$, the integration by parts of
the operator
$\Delta +\partial_t$, which is acting on $u^2$ over the other terms in the corresponding
integral over $\Rn\times [0,4)$ and using the Cauchy-Schwarz inequality to handle the
cross term, \eqref{E:acotaciones} and \eqref{E:carleman2}.  
\end{proof}
\begin{lemma}\label{L:decay3}
Given $\theta\in (0,1]$, there are $N=N(n,\theta)\ge 1$ and  $\rho=\rho(n,\theta)$ in $(0,1]$
such that the following holds:

If $u$ satisfies
$|\Delta u+\partial_tu |\le  R^2|u|+R|\nabla u|$,
$\|u\|_{\infty}\le 1$,
$\|\nabla u\|_\infty\le R$ in $B_4\times [0,\frac 1{R^2}]$ and 
$\|u(\,\cdot\,,0)\|_{L^2\left(B_\rho\right)}\ge \theta R^{-n/2}$. Then, 
\begin{itemize}
\item $\|u(\,\cdot\,,0)e^{-R^2|x|^2/\epsilon}\|_{L^2\left(B_{4}\right)}\ge 
e^{-NR^2\log{\left(\frac 1{\epsilon}\right)}}$, when $0<\epsilon\le \tfrac 1{3N}\ $, 
$R\ge N$.
\item $\|u(\,\cdot\,,0)\|_{L^2\left(B_{r}\right)}\ge 
e^{-NR^2\log{\left(\frac Nr\right)}}$, when $0<r\le \frac 12\ $, $R\ge N$.
\end{itemize}
\end{lemma}
\begin{proof} 
Take as $f$ in Lemma
\ref{L:carleman1} the function,
$f=u \varphi(x)\psi(t)$, where
$\varphi\in C_0^\infty(B_4)$,
$0\le \varphi\le 1$, $\varphi=1$ in $B_3$ and $\varphi=0$ outside $B_{\frac 72}$, $\psi=1$ when
$0<t\le\frac 1\alpha$ and
$\psi=0$ when $t\ge\frac 2\alpha$.  Then,
\begin{equation}\label{E:error}
|\Delta f+\partial_tf |\le R^2|f|+R|\nabla f|+N\left(\al+R\right)\chi_{B_4\times [0,\frac 2\al]\setminus
B_{3}\times [0,\frac 1\al]}\ . 
\end{equation}
The facts that  $t/N\le \sigma (t)\le t$ on $[0,\frac 6\al]$, that $\s^{1-\al}_aG_a\le
N^{\al+\frac n2}\al^{\al+\frac n2-1}$ in the region
$B_4\times [0,\frac 2\al]\setminus B_{3}\times [0,\frac 1\al]$, \eqref{E:error} and standard arguments
with Carleman inequalities imply the estimate
\begin{multline}\label{E:applicarleman}
\al^2\int_0^{\frac 1\al}\int_{B_2}\left(t+a\right)^{-\al}u^2e^{-|x|^2/4(t+a)}\,dX\le
N^{\al}\al^{\al+1}\\+N^\al\s (a)^{-\al}\left [-(a/N)\int |\nabla
f(x,0)|^2e^{-|x|^2/4a}\,dx+N\al\int f ^2(x,0)e^{-|x|^2/4a}\,dx\right]\ ,
\end{multline}
when
$\al\ge NR^2$ and $0<a\le \frac 1\al\ $.\par
For $\rho$ in $(0,1]$, which will be chosen later and Lemma \ref{L:decay}, we know that
\[
\sqrt N\|u(\,\cdot\,,t)\|_{L^2\left(B_{2\rho}\right)}\ge 
R^{-n/2},\ \text{when}\  0<t\le 1/R^2\  \text{and}\  R> N/\rho\ .
\]
This and the
conditions, $0<a\le
\frac{\rho^2}{2\al}$ and $\al\ge NR^2$, imply that the left hand side of \eqref{E:applicarleman} is
bounded from below by 
\begin{equation}\label{E:applicarleman22}
\al^2\int^{\frac
{\rho^2}{\al}-a}_{\frac
{\rho^2}{2\al}-a}\int_{B_{2\rho}}\left(t+a\right)^{-\al}e^{-\frac{\rho^2}{(t+a)}}u^2\,dX
\ge\frac{\al^{\al+1}\rho^2}{2NR^n}
\left(\frac
1{\rho e}\right)^{2\al}\ .  
\end{equation}
Inequalities \eqref{E:applicarleman} and \eqref {E:applicarleman22} show, that  to make
sure that the left hand side of
\eqref{E:applicarleman} is larger than  four times the  first term on right hand side of
\eqref{E:applicarleman}, when $\al\ge NR^2$ and $0<a\le
\frac{\rho^2}{2\al}$, it suffices to know that
\begin{equation}\label{E:suffices}
\left(\frac 1{\rho e}\right)^{2\al}\ge 8N^{\al+1} R^n/\rho^2\ .   
\end{equation}
Choose then $\rho$ as the solution of the equation $\frac 1{\rho e}=\sqrt {8N}$.
Then,
\eqref{E:suffices}
 holds when $8^{\al-1}\ge NR^n/\rho^2$. Thus, there are fixed constants,
$\rho=\rho(n,\theta)$ in $(0,1]$ and $N=N(n,\theta)\ge 1$ such that, under the conditions
in Lemma \ref{L:decay3}, we have  
\begin{multline*}\label{E:resultado}
\tfrac 12\int_0^{\frac
1\al}\int_{B_2}\left(t+a\right)^{-\al}u^2e^{-|x|^2/4(t+a)}\,dX+N^{\al}\al^{\al+1}\\\le N^{\al}\s
(a)^{-\al}\left [-(a/N)\int |\nabla f(x,0)|^2e^{-|x|^2/4a}\,dx+N\al\int f
^2(x,0)e^{-|x|^2/4a}\,dx\right]\ ,
\end{multline*}
when $R\ge N$, $\al\ge NR^2$ and $0<a\le \frac{\rho^2}{12\al}\ $. In particular, there is
$N=N(n,\theta)$ such that
\begin{equation}\label{E:resultado}
N^{-\al}\al^{\al+1}a^\al\le -(a/N)\int |\nabla f(x,0)|^2e^{-|x|^2/4a}\,dx+N\al\int f
^2(x,0)e^{-|x|^2/4a}\,dx\ ,
\end{equation}
when $R\ge N$, $\al=NR^2$ and $0<a\le \frac 1{12NR^2}\ $.

Recalling the definition of
$f$, choose $a=\frac \epsilon{8R^2}$ in \eqref{E:resultado}. It implies the inequality
\begin{equation*}\
e^{-2NR^2\log{\left(\frac 1{\epsilon}\right)}}\le
\int_{B_4}u^2(x,0)e^{-2R^2|x|^2/\epsilon}\,dx\ ,\ \text{when}\ 0<\epsilon\le \tfrac 2{3N}\ ,
\ R\ge N 
\end{equation*}
and proves the first claim in Lemma \ref{L:decay3}.
The inequality \eqref{E:resultado} also implies the bound 
\begin{multline*}
2a\int |\nabla f(x,0)|^2e^{-|x|^2/4a}\,dx+\tfrac n2\int
f^2(x,0)e^{-|x|^2/4a}\,dx\\\le NR^2\int
f^2(x,0)e^{-|x|^2/4a}\,dx\quad ,\quad \text{when}\  0< a\le \tfrac 1{12NR^2}\ ,\ 
R\ge N\ .
\end{multline*}
From Lemma \ref{L:implidoub} with $h=f(\ \cdot,0)$ and the above estimate, we obtain 
\begin{equation}\label{E:doubling}
\int_{B_{2r}}u^2(x,0)\,dx\le e^{NR^2}\int_{B_r}u^2(x,0)\,dx\ , \text{when}\ 0<r\le
1/2\ ,\ R\ge N\ .
\end{equation}
For these values of $r$,
choose
$k\ge 2$ such that,
$2^{-k}< r\le 2^{-k+1}$ and iterate \eqref{E:doubling} when $r=2^{-j}$, $j=0, \dots ,k-1$. It gives 
\begin{equation*}
\int_{B_1}u^2(x,0)\,dx\le e^{2NR^2\log{\left(1/r\right)}}\int_{B_r}u^2(x,0)\,dx\ ,
\text{when}\ 0<r\le 1/2\ ,\ R\ge N\ ,
\end{equation*}
which proves the second claim.
\end{proof}
\begin{proof}[Proof of Theorem 1] Without loss of generality we may assume that $u$
satisfies
\begin{equation*}
|\Delta u+\partial_tu|\le |u|+|\nabla u|\ \text{and}\ \  |u|\le
1\ \text{in}\ \ \R^n\times [0,4]\ .
\end{equation*}
Choose $\theta$ in $(0,1]$ such that, $\theta\le
\|u(\,\cdot\,,0)\|_{L^2\left(B_1\right)}$. If 
$\rho$ is the constant associated to $\theta$ in Lemma
\ref{L:decay3} , define $u_R(x,t)=u(Rx+y,R^2t)$, when   
$R\rho=2|y|$ is large, $y\in\Rn$. Then,  
\[ R^{n/2}\|u_R(\,\cdot\,,0)\|_{L^2\left(B_{\rho}\right)}
=\|u(\,\cdot\,,0)\|_{L^2\left(B_{2|y|}(y)\right)}\ge
\|u(\,\cdot\,,0)\|_{L^2\left(B_1\right)}\ge \theta
\] 
and the standard interior estimates for solutions to parabolic equations \cite{ls} show that
$u_R$ satisfies the conditions in Lemma
\ref{L:decay3}. The first claim in Lemma \ref{L:decay3} applied to $u_R$ and the change of
variables, $Rx+y=z$, give that for
$\epsilon$ sufficiently small

\begin{multline}\label{E:final}
R^{n/2}e^{-NR^2\log{\left(\frac
1{\epsilon}\right)}}\le
\|u(\,\cdot\,,0)e^{-|x-y|^2/\epsilon}
\|_{L^2\left(B_{4|y|/\rho}(y)\right)}\\
\le \|u(\,\cdot\,,0)
\|_{L^2\left(B_{|y|/2}(y)\right)}+R^{n/2}e^{-R^2/8\epsilon}\ ,\text{ when}\ 0<\epsilon\le \tfrac
1{N}\ , \ R\ge N\ ,
\end{multline}
 and choosing $\epsilon$ small in \eqref{E:final} implies the first
inequality in Theorem \ref{T:1}.

The second claim in Lemma \ref{L:decay3} applied to $u_R$ with $r=1/R$ and the same change
of variables, give
$$R^{-n/2}\|u(\,\cdot\,,0)\|_{L^2\left(B_{1}(y)\right)}=\|u_R(\,\cdot\,,0)\|_{L^2\left(
B_{
1/R}\right)}\ge  e^{-NR^2\log{\left(NR\right)}}\ ,$$
which proves the second inequality in Theorem \ref{T:1}. 

What has been proved so far shows that the condition
\[|u(x,0)|\le C_ke^{-k|x|^2}\ \text{for all}\ 
x\in\Rn\ \text{and}\  k\ge 1\] can only hold when
$\|u(\,\cdot\,,0)\|_{L^2\left(B_{1}\right)}$ vanishes. The results in \cite{av} or
\cite[Theorem 3]{f} prove that the latter is only possible, when
$u(\,\cdot\,,0)\equiv 0$. Then, standard backward uniqueness arguments for
parabolic equations imply,
$u\equiv 0$ in
$\Rn\times [0,4]$, when
$u\in L^\infty(0,4\ ;L^2(\Rn))$ \cite{evans}. If one wants to relax the latter
condition and to allow
$u$ to grow as a quadratic exponential at infinity in the layer $\Rn\times [0,4]$, the fact
that $u\equiv 0$ in
$\Rn\times [0,4]$ follows from the arguments in \cite[Theorem 3]{f} or  the
Carleman inequality \eqref{E:primera desigualdad}. 
\end{proof} 
\section{Second Proof of Theorem \ref{T:1}}\label{S:3}
The second proof of Theorem \ref{T:1} is based in Lemmas \ref{L:freq} and \ref{L:ineq}. 
\begin{lemma}\label{L:freq}
Given $a>0$ and $f\in W^{2,\infty}(\RM)$, set
\[H_a(t)=\int_{\Rn}f^2G_a\,dx\ ,\ D_a(t)=\int_{\Rn}|\nabla f|^2G_a\,dx
\ \text{and}\  N_a(t)=\frac{2(t+a)D_a(t)}{H_a(t)}\ ,
\]
where $G_a(x,t)=(t+a)^{-n/2}e^{-|x|^2/4(t+a)}$.
Then, 
\begin{equation*}
\dot N_a(t)\ge -\frac{(t+a)}{H_a(t)}\int (\Delta
f+\partial_tf)^2G_a\, dx\quad .
\end{equation*}
\end{lemma}
The monotonicity results implied by this Lemma (e.g. $N_a(t)$ is nondecreasing when
$f$ is a backward caloric in $R^{n+1}_+$) are within the category of what in the literature
have been called  frequency function arguments. The frequency function here is
$N_a(t)$. This frequency function seems to have first 
appeared or been used in the context of unique continuation for parabolic equations in \cite{p},
when
$a=0$ and in
\cite{f}, when
$a>0$. Related results, though with perhaps different
purposes, appeared in \cite{ha} and \cite{ha2}. 

The next proof of Lemma \ref{L:freq} comes from \cite[Lemma 2]{efv}.
\begin{proof}
The identities $\partial_tG_a-\Delta G_a=0$, $\nabla G_a=-\frac x{2\left(t+a\right)}G_a$,
$\Delta=\text{div}\left(\nabla\ \right)$ and integration by parts imply the following
identities
\begin{equation}\label{E:deriv}
\dot H_a(t)=2\int f(\Delta f+\partial_tf)G_adx+2D_a(t)\ , 
\end{equation}
\begin{align*}\dot
H_a(t)&=2\int f\left(\partial_tf+\tfrac x{2(t+a)}\cdot
\nabla f-\tfrac 12\left(\Delta f+\partial_tf\right)\right)G_a\,dx+\int f\left(\Delta f+\partial_t
f\right)G_a\,dx\ ,\\ D_a(t)&=\int f\left(\partial_tf+\tfrac {x}{2(t+a)}\cdot \nabla
f-\tfrac 12 \left(\Delta f+\partial_tf\right)\right)G_a\,dx-\tfrac 12\int f\left(\Delta f+\partial_t
f\right)G_a\,dx
\end{align*}
and
\begin{align}\label{E:mult}
\dot H_a(t)D_a(t)&=2\left(\int f\left(\partial_tf+\tfrac x{2(t+a)}\cdot \nabla f-\tfrac 12\left(\Delta
f+\partial_tf\right)\right)G_a\,dx\right)^2\\&-\tfrac 12\left(\int f\left(\Delta f+\partial_t
f\right)G_a\,dx\right)^2\ .\notag
\end{align}
The Rellich-N\v{e}cas identity with vector field $\nabla G_a$ 
\begin{align*}
\text{div}(\nabla G_a|\nabla f|^2)&-2\text{div}((\nabla f\cdot\nabla G_a)\nabla
f)\\&=\Delta G_a|\nabla f|^2-2D^2G_a\nabla f\cdot\nabla f-2\nabla f\cdot\nabla G_a\Delta f
\end{align*}
and integration by parts give
\begin{align}\label{E:Rellich}
\int\Delta G_a&|\nabla f|^2\,dx=2\int D^2G_a\nabla f\cdot\nabla f\,dx+2\int\nabla
f\cdot\nabla G_a\Delta f\,dx\\&=2\int \left(\tfrac x{2(t+a)}\cdot\nabla
f\right)^2G_a\,dx-2\int\tfrac x{2(t+a)}\cdot\nabla f\Delta fG_a\,dx-D_a(t)/(t+a)\ .\notag
\end{align}
Again, the fact that $G_a$ is a caloric function, integration by parts, \eqref{E:Rellich} and the completion of
the square of $\partial_tf+\tfrac {x}{2(t+a)}\cdot
\nabla f-\tfrac 12\left(\Delta f+\partial_tf\right)$ yields the formula
\begin{align}\label{E:derv2}
\dot D_a(t) &=2\int\left(\partial_tf+\tfrac
{x}{2(t+a)}\cdot
\nabla f-\tfrac 12\left(\Delta f +\partial_tf\right)\right)^2G_a\,dx\\&-\tfrac 12\int\left(\Delta
f+\partial_tf\right)^2G_a\,dx-D_a(t)/(t+a)\ .\notag
\end{align}
Then, from \eqref{E:mult},\eqref{E:derv2} and the quotient rule
\begin{align}\label{derv5}
\dot N_a(t)&=\frac {4(t+a)}{H_a(t)^2}\Biggl\{\int\left(\partial_tf+\tfrac
{x}{2(t+a)}\cdot\nabla f-\tfrac 12 \left(\Delta
f+\partial_tf\right)\right)^2G_a\,dxH_a(t)\\ &-\left(\int f\left(\partial_tf+\tfrac
{x}{2(t+a)}\cdot \nabla f-\tfrac 12\left(\Delta
f+\partial_tf\right)\right)G_a\,dx\right)^2\notag\\& +\tfrac{1}{4}\left(\int f\left(\Delta
f+\partial_tf\right)G_a\,dx\right)^2-\tfrac{1}4\int\left(\Delta
f+\partial_tf\right)^2G_a\,dxH_a(t)\Biggr\}\notag
\end{align}
and Lemma \ref{L:freq} follows from \eqref{derv5}, the Cauchy-Schwarz inequality and the
positiveness of the third term on the right hand side of \eqref{derv5}.
\end{proof}
The application of the Lemmas \ref{L:freq}  and \ref{L:ineq} to the proof of Theorem \ref{T:1} is the
following: 
\begin{lemma}\label{L:decay2}
Assume that $u$ satisfies, $\left |\Delta u+\partial_tu\right |\le  R^2|u|+R|\nabla u|$,
$\|u\|_{\infty}\le 1$,
$\|\nabla u\|_\infty\le R$ in $B_4\times [0,\frac 1{R^2}]$ and 
$\|u(\,\cdot\,,0)\|_{L^2(B_1)}\,dx\ge \theta R^{-n/2}$
for some $\theta\in (0,1]$. Then, there is $N=N(n,\theta)$ such that
\begin{itemize}
\item $\|u(\,\cdot\,,0)e^{-R^2|x|^2/\epsilon}\|_{L^2\left(B_{4}\right)}\ge 
e^{-NR^2\log{\left(\frac 1{\epsilon}\right)}}$, when $0<\epsilon\le \tfrac 1{3N}\ $, 
$R\ge N$.
\item $\|u(\,\cdot\,,0)\|_{L^2\left(B_{r}\right)}\ge 
e^{-NR^2\log{\left(\frac Nr\right)}}$, when $0<r\le \frac 12\ $, $R\ge N$.
\end{itemize}
\end{lemma}
\begin{proof} The Lemma \ref{L:decay} with $\rho=1$ gives
\begin{equation}\label{E:inferior}
\sqrt{N}\|u(\,\cdot\,,t)\|_{L^2\left(B_{2}\right)}\ge R^{-n/2}\ ,\ \text{when}\ 
0<t\le \tfrac 1{R^2}\ \text{and}\ R\ge N\ .
\end{equation}
Set 
$f=u\varphi$ in Lemma
\ref{L:freq}, where
$\varphi\in C_0^\infty(B_4)$,
$0\le\varphi\le 1$, $\varphi=1$ in $B_3$ and $\varphi=0$ outside $B_{\frac 72}$. From
\eqref{E:inferior},
\begin{equation}\label{E:necea}
H_a(t)\ge N^{-1}R^{-n}
(t+a)^{-n/2}e^{-1/(t+a)}\ ,\ \text{when}\ t+a\le\tfrac 1{R^2}\ ,\ R\ge N\ .
\end{equation}
and
\begin{equation}\label{E:heat4}
|\Delta f+\partial_tf|\le R^2|f|+R|\nabla f|+NR\chi_{B_4\setminus B_{3}}\ ,\ \text{in}\
B_4\times [0,\tfrac 1{R^2}]\ .
\end{equation}
From  Lemma \ref{L:freq}, \eqref{E:heat4} and \eqref{E:necea}, we
have
\begin{equation}\label{E:neceb}
\dot N_a(t)\ge -NR^2-NR^2N_a(t)\ ,\ \text{when}\ t+a\le\tfrac 1{R^2}\ ,\ R\ge N\ .
\end{equation} 
Thus
\begin{equation}\label{E:incres}
e^{NR^2t}N_a(t)+e^{NR^2t}\ \text{is nondecreasing,
when}\ t+a\le \tfrac 1{R^2}\ \text{and}\ R\ge N\ .
\end{equation}
The multiplication of the identity
\begin{equation*}
\dot{H}_a(t)=2\int f(\Delta f+\partial_tf)G_adx+2D_a(t)
\end{equation*}
by $(t+a)/H_a(t)$, \eqref{E:heat4} and \eqref{E:necea}, imply that for some  $N>0$,
\begin{equation}\label{E:relacion}
N_a(t)\le N\left[1+(t+a)\partial_t\log {H_a(t)}\right]\ ,\ \text{when}\ \ 0\le t+a\le
\tfrac 1{R^2}\ .
\end{equation}
Set $\beta =\frac 1{R^2}\ $. Then, from \eqref{E:relacion}, \eqref{E:incres} and
\eqref{E:necea} 
\begin{equation*}
\begin{split}
N_a(0)&\lesssim N_a(\beta/4)+1\lesssim 1+\int_{\beta /4}^{\beta
/2}\frac{N_a(t)}{(t+a)}\,dt\lesssim 1+\int_{\beta /4}^{\beta
/2}\partial_t\log{H_a(t)}\,dt\\&=1+\log{\left(\tfrac{H_a(\beta/2)}{H_a(\beta/4)}\right)}
\le NR^2\ ,\ \text{when}\ a\le\tfrac \beta{12}\ .
\end{split}
\end{equation*}
In particular,
\begin{equation}\label{E:finale}
\begin{split}
2a&\int|\nabla f(x,0)|^2e^{-|x|^2/4a}\,dx+\tfrac n2\int f^2(x,0)e^{-|x|^2/4a}\,dx\\ &\le
 NR^2\int f^2(x,0)e^{-|x|^2/4a}\,dx\ , \
\text{when} \ 0<a\le\tfrac {1}{12NR^2}\ ,\  R\ge N\ .
\end{split}
\end{equation}
Now, Lemma \ref{L:ineq} and \eqref{E:finale}
give
\[\partial_a\log{\left(\int f^2(x,0)e^{-|x|^2/4a}\,dx\right)}\le \tfrac{NR^2}{a}\ , \
\text{when} \ 0<a\le\tfrac {1}{12NR^2}\ ,\  R\ge N\ , 
\]
and the integration of this inequality over
$[\tfrac{\epsilon}{4R^2},\tfrac 1{12NR^2}]$ implies the first claim in Lemma \ref{L:decay2}.
The second claim is derived from
\eqref{E:finale}, as in Lemma \ref{L:decay3}.  
 \end{proof}
\begin{proof}[Proof of Theorem 1] Proceeding as before, we may  assume that 
$u$ satisfies
\begin{equation*}
|\Delta u+\partial_tu|\le |u|+|\nabla u|\ \text{and}\ \  |u|\le
1\ \text{in}\ \ \R^n\times [0,4]\ .
\end{equation*}
Choose then $\theta$ in $(0,1]$ such that, $\theta\le
\|u(\,\cdot\,,0)\|_{L^2\left(B_1\right)}$ and set $u_R(x,t)=u(Rx+y,R^2t)$, when   
$R=2|y|$ is large, $y\in\Rn$. Then,  
\[ R^{n/2}\|u_R(\,\cdot\,,0)\|_{L^2\left(B_{1}\right)}
=\|u(\,\cdot\,,0)\|_{L^2\left(B_{2|y|}(y)\right)}\ge
\|u(\,\cdot\,,0)\|_{L^2\left(B_1\right)}\ge \theta
\] 
and the arguments proceed as in the first proof. 
\end{proof}
\section{Proof of Theorem \ref{T:2}}\label{S:4}
Theorem \ref{T:2} follows with similar arguments. Here is an outline of its proof with a
frequency function type argument .
\begin{proof}[Proof of Theorem 2] As before and without loss of generality we may assume that
\begin{equation*}
|\Delta u+\partial_tu|\le |u|+|\nabla u|\ \text{and}\ \  |u|\le
1\ \text{in}\ \ \R^n_+\times [0,4]\ .
\end{equation*}
Choose then $\theta$ in $(0,1]$ such that, $\theta\le
\|u(\,\cdot\,,0)\|_{L^2\left(B_1(4e_n)\right)}$. The argument in the proof of Lemma
\ref{L:decay} (See \cite[Lemmas 1]{efv}) is easily adapted to
show that there is $N=N(n,\theta)>0$ such that
\begin{equation}\label{E:por abajo}
\sqrt N\|u(\,\cdot\,,t)\|_{L^2\left(B_{2}(4e_n)\right)}\ge 1\ ,\ \text{when}\ 
0<t\le \tfrac 1{N}\ .
\end{equation}
Set 
$v(x,t)=u(yx+ye_n,y^2t)$, when $y>8$. The fact that the ball of radius $2/y$
and centered at $(4/y)e_n-e_n$ is contained in
$B_{1-2/y}$, the change of variables, $z=yx+ye_n$, and \eqref{E:por abajo}, imply that
\begin{equation}\label{E:inferior22}
\sqrt N\|v(\,\cdot\,,t)\|_{L^2\left(B_{1-2/y}\right)}\ge
y^{-n/2}\ ,\ \text{when}\ 0\le t\le 1/(Ny^2)\ .
\end{equation}
Set 
$f=v\varphi$ in Lemma
\ref{L:freq}, where
$\varphi\in C_0^\infty(B_1)$,
$0\le\varphi\le 1$, $\varphi=1$ in $B_{1-1/y}$ and $\varphi=0$ outside $B_{1-1/(2y)}$. From
\eqref{E:inferior22}
\begin{equation}\label{E:necea2}
H_a(t)\ge N^{-1}y^{-n}
(t+a)^{-n/2}e^{-\left(1-\tfrac 2y\right)^2/4(t+a)}\ ,\ \text{when}\ t+a\le\tfrac
1{Ny^2}\ .
\end{equation}
and
\begin{equation}\label{E:heat22}
|\Delta f+\partial_tf|\le y^2|f|+y|\nabla f|+Ny^2\chi_{B_1\setminus B_{1-1/y}} \ ,\ 
\text{in}\ B_1\times [0,4/y^2]\ .
\end{equation}
The calculations, which were carried out out in the second proof of Theorem
\ref{T:1} but replacing \eqref{E:necea} and \eqref{E:heat4} by
\eqref{E:necea2} and
\eqref{E:heat22} respectively, imply the inequality
\begin{equation}\label{E:finale22}
\begin{split}
2a\int |\nabla f(x,0)|^2e^{-|x|^2/4a}\,dx+\tfrac n2\int
f^2(x,0)e^{-|x|^2/4a}\,dx\\\le Ny^2\int
f^2(x,0)e^{-|x|^2/4a}\,dx\ , \ \text{when} \ 0<a\le\tfrac {1}{12Ny^2}\ , 
\end{split}
\end{equation}
which as seen before, implies the first part of Theorem
\ref{T:2}. The second claim follows from the first, the results in \cite{av} or
\cite[Theorem 3]{f} and the qualitative result in \cite{ess2}, which was stated in the Introduction after Theorem
\ref{T:1}.   
\end{proof}

\newpage

\end{document}